\documentclass{article}

\usepackage{amsmath,amssymb,amsthm}
\usepackage[utf8]{inputenc}
\usepackage{color,graphicx}

\newcommand{\scp}[2]{{\left\langle {#1}\, , \, {#2}\right\rangle}}

\newcommand{\mR}{{\mathbb R}}

\newcommand{\prt}{\partial}
\newcommand{\sub}{\subset}

\renewcommand{\phi}{\varphi}
\newcommand{\Om}{\Omega}
\newcommand{\ga}{\gamma}
\renewcommand{\th}{\theta}
\newcommand{\la}{\lambda}
\newcommand{\al}{\alpha}
\newtheorem{proposition}{Proposition}

\graphicspath{{./normalCortex/}}

\title{Normal and Equivolumetric  Coordinate Systems for Cortical Areas}
\author{Laurent Younes$^1$ \and Kwame S. Kutten$^2$ \and J. Tilak Ratnanather$^2$}
\date{\bigskip
$^1$Department of Applied Mathematics and Statistics, Johns Hopkins University\\
$^2$Department of Biomedical Engineering, Johns Hopkins University\\
\bigskip
\today}

\begin{document}
\maketitle

\begin{abstract}
 We describe coordinate systems adapted for the space between two surfaces, such as those delineating the highly folded cortex in mammalian brains. These systems are estimated in order to satisfy geometric priors, including streamline normality or equivolumetric conditions on layers. We give a precise mathematical formulation of these problems, and present numerical simulations based on diffeomorphic registration methods, comparing them with recent approaches.
\end{abstract}

\section{Introduction}
Anatomical regions of interest extracted from 3D biomedical imaging data often appear as volumes separated by ``upper'' and ``lower'' surfaces. These include, in particular, the highly folded cortical regions or areas of the mammalian brain. This paper focuses on methods that parametrize such volumes, using coordinate systems that are naturally aligned with the encompassing surfaces. Our contributions are, on one hand, to formalize a notion of laminar coordinate systems and, on the other hand, discuss within this formalism the concept of equivolumetric coordinates \cite{Bok1929,Bok1959}, providing an interpretation of two similar methods \cite{waehnert2014anatomically,leprince2015combined} and introducing a new one based on diffeomorphic registration methods.

\begin{figure}
    \centering
\includegraphics[width=0.4\textwidth]{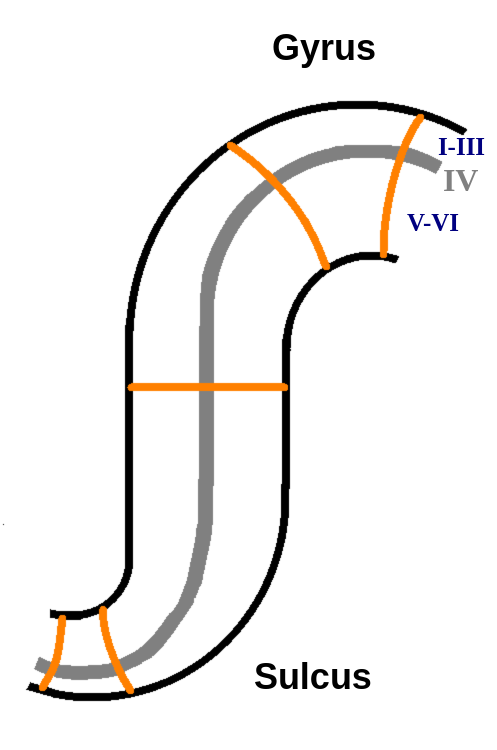}
   \caption{Idealisation of Bok's equivolumetric hypothesis for the layers in the folded cortex showing that the deep layers (V-VI) are thin in the sulcus and thick in the gyrus respectively characterized by regions of positive and negative inward curvature and conversely for the upper layers (I-III)}
\label{fig:bok}
\end{figure}

The elegant laminar structure of a typical cortical area \cite{Dahnke2018,Wagstyl2018} is summarized as follows. The folding of the area serves to maximize its surface area in a confined cranial space. The neural tissue (grey matter) within the area contains mostly neuronal cell bodies and unmyelinated fibers. Cortical areas (which number in the hundreds in the human brain) are connected via white matter containing axonal, usually myelinated, fibers. Each cortical area is composed of fundamental units called cortical columns that traverse vertically from the white matter to the surface just below the pial matter. Finally, the cortical area is composed of six layers which are stacked horizontally on top of each other. Bok \cite{Bok1929,Bok1959} observed that to maintain the laminar structure in highly folded regions that thin layers in one part became thicker in another part. This observation led to the hypothesis that the cortex satisfies an equivolumetric property,  illustrated in Figure \ref{fig:bok}, which provides the theoretical motivation of this paper.

The rest of the paper is as follows. Section \ref{sec:laminar} introduces formal definitions of laminar coordinates and describes methods for constructing them. Section \ref{sec:bok} describes how Bok's equivolumetric hypothesis can be implemented. Section \ref{sec:results} presents
numerical simulations.

\section{Laminar Coordinate Systems and Thickness}
\label{sec:laminar}
\subsection{Notation}
\begin{figure}
    \centering
    \includegraphics[trim=2cm 12cm 2cm 3cm, clip,width=0.9\textwidth]{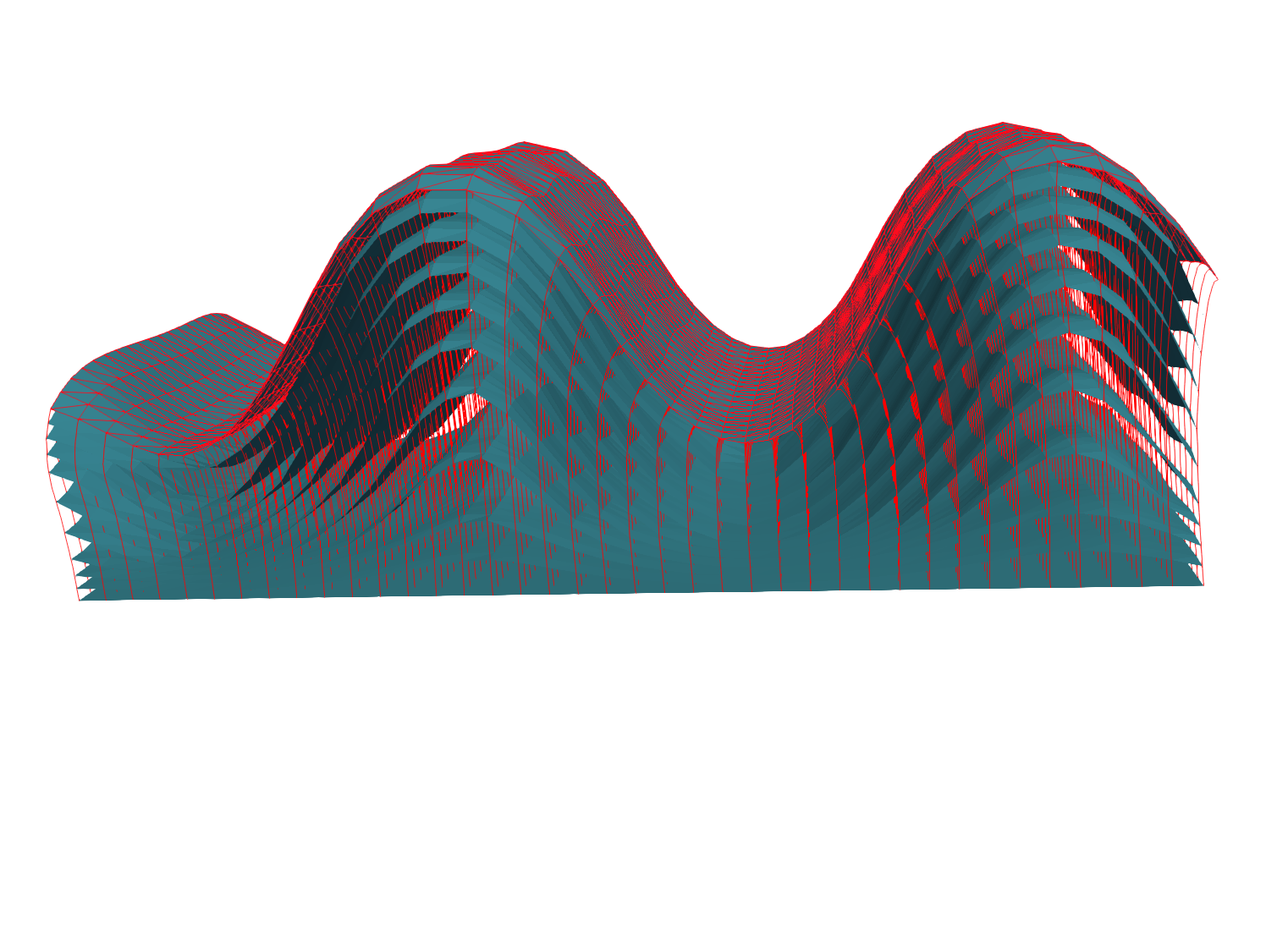}
    \caption{A laminar coordinate system between two surfaces, represented by the red scaffold.}
    \label{fig:waves}
\end{figure}

We introduce some mathematical notation in order to describe 3D coordinate systems parametrizing an open space between two surfaces. We will call ``laminar coordinate system'' (see Figure \ref{fig:waves} for an example) a special case of a {\em foliation} of this open set, with two special leaves provided by the two surfaces.

{}

More precisely, let $S_0$ and $S_1$ be 2D submanifolds (surfaces) in $\mR^3$, such that $S_0\cap S_1 = \emptyset$.  Assume that these surfaces are bounded, and that they are either compact  (e.g., spheres), or surfaces with boundaries (e.g., disks). A laminar coordinate system between $S_0$ and $S_1$ is a $C^1$ embedding $\psi: [0,1] \times S_0 \to \mR^3$ such that
$\psi(0,x) = x$ for all $x\in S_0$ and $\psi(1, \cdot)$ maps $S_0$ onto $S_1$.
The embedding assumption requires that (i) $\psi$ in one-to-one,  (ii) $\partial_x\psi(t,x)$ (the differential of $\psi(t, \cdot)$ at $x$) can be extended by continuity to the whole set $[0,1] \times S_0$ and is an invertible linear mapping, and (iii) the inverse mapping defined on $\bar\Omega = \psi([0,1]\times S_0)$ is continuous.

The surfaces $S_t = \psi(t, S_0)$ are the leaves of the foliation, and will be referred to as ``layers'' while the curves $\ga_x(t) = \psi(t, x)$, $t\in [0,1]$ will be referred to as ``streamlines.''  These layers and streamlines are shown as blue surfaces and vertical red lines respectively in Figure \ref{fig:waves}. The set $\bar\Omega = \psi([0,1]\times S_0)$ is the closure of an open subset $\Omega$ of $\mR^3$ defining the space between the two surfaces.

In addition, a coordinate system defined as such provides a definition of ``thickness'' of $\Omega$ at $x\in S_0$, given by the length of the streamline starting at $x$, namely,
\[
\th(x)= \int_0^1 |\prt_t \psi(t,x)|\, dt.
\]

\subsection{Coordinate Systems}
There are several methods for building coordinate systems. Three approaches are described as follows.
\bigskip

\subsubsection{Level Set Methods}
Assume that $\Om\sub \mR^3$ and that its boundary has two connected components, providing $S_0$ and $S_1$, which must therefore be closed surfaces. The typical application is when $\Omega$ is the region between two non-intersecting nested sphere-like surfaces.

Assume that one is given a $C^2$ submersion $F: \bar \Om\to [0,1]$ such that $S_0 = F^{-1}(0)$ and $S_1 = F^{-1}(1)$. (This assumes that $F$ is onto and that $\nabla F(x)\neq 0$ for all $x\in \bar\Om$.) Define the vector field
\[
v(x) = \nabla F/|\nabla F|^2,
\]
and let $\psi(t,x)$, $t\in [0,1]$, $x\in S_0$ satisfy
\[
\prt_t \psi(t, x) = v(\psi(t,x))
\]
with $\psi(0,x) = x$, so that $\psi$ is the restriction to $S_0$ of the flow associated with the vector field $v$. Then, one has
\[
\prt_t F(\psi(t,x)) = \nabla F(\psi(t,x))^T v(\psi(t,x)) = 1
\]
which implies that $F(\psi(t,x)) = t$ for all $x \in S_0$ and $t\in [0,1]$. The function $\psi: [0,1]\times S_0 \to \bar\Om$ provides a laminar coordinate system of $\bar\Om$ as defined above. The thickness along streamlines is in particular defined by
\[
\th(x) = \int_0^1 |\nabla F(\psi(t,x))|^{-1}\, dx.
\]
Notice that, since the layers $S_t$ coincide with the level sets of $F$, the vector field $v$ (and therefore the streamlines) are necessarily perpendicular to them.
\medskip

Given $S_0$ and $S_1$, one therefore needs to define a suitable function $F$. Jones et al. \cite{jones2000three} defined $F$ as the solution of the Laplace equation $\Delta F = 0$ with boundary conditions $F=0$ on $S_0$ and $F=1$ on $S_1$. Similarly, Waehnert et al. \cite{waehnert2014anatomically} proposed a construction starting with level-set representations of $S_0$ and $S_1$, represented by two functions $U_0$ and $U_1$ such that $S_i = U_i^{-1}(0)$, $i=0,1$. This leads to the definition
\begin{equation}
    \label{eq:waehnert1}
F(t, \cdot) = t + \mathcal S((1-\rho)U_0 + \rho U_1)
\end{equation}
where $\mathcal S$ is a smoothing operator based on a topology-preserving mean-curvature motion equation. More precisely, given a function $U_{\mathrm{targ}}$, the function $\mathcal S U_{\mathrm{targ}}$ is obtained as the limit when time $s\to\infty$ of the solution $U(s, \cdot)$ of
\[
\prt_s U + (U-U_{\mathrm{targ}}) |\nabla U| = \epsilon |\nabla U|\mathrm{div}\left(\frac{\nabla U}{|\nabla U|}\right)\,.
\]
\bigskip

\subsubsection{Diffeomorphic Volume Mapping}
Das et al. \cite{das2009registration} assumed that the surfaces $S_0$ and $S_1$ are represented as boundaries of two open set $\Om_0$ and $\Om_1$ (so that $S_0 = \prt\Om_0$ and $S_1=\prt \Om_1$)  with $\bar\Om_0 \sub \Om_1$. Then a variant \cite{Avants2004} of the large deformation diffeomorphic metric mapping (LDDMM) algorithm \cite{beg2005computing} was used to estimate a flow of diffeomorphisms $\phi(t, \cdot): \mR^3\to \mR^3$, $t\in [0,1]$ such that $\phi(1, \Om_0) \simeq \Om_1$.
The function  $\psi$ can then be defined as the restriction of $\phi$ to $[0,1] \times S_0$.
\bigskip

\subsubsection{Surface Mapping with Normality Constraints}
The foregoing methods are volumetric: they work on 3D regions rather than directly on surfaces, which are represented either as level sets or as boundaries of open sets. They are, in particular, not well adapted to analyze regions delimited by open surfaces. So following \cite{ratnanather20193D}, a version of LDDMM adapted to surface mapping was combined with normality constraints in order to estimate laminar coordinate systems. This is now described in some detail.
\medskip

Assume that $S_0$ is parametrized over a bounded open subset $M\subset \mR^2$ (or more generally over a 2D manifold $M$ with or without boundary) in the form $S_0 = q_0(M)$, where $q_0$ is an embedding of $M$ into $\mR^3$.  The LDDMM surface registration algorithm solves an optimal control problem minimizing, over all time-dependent vector fields in a reproducing kernel Hilbert space $V$,
 \[
 \int_0^1 \|v(t)\|_V^2\, dt + D(q(1, M), S_1)
 \]
subject to $q(0) = q_0$ and $\prt_t q(t) = v(t) \circ q(t)$. Here,
$D$ is a reparametrization-invariant discrepancy measure between (unparametrized) surfaces. Several versions of this cost function have been introduced, based on representation of surfaces as currents \cite{vaillant2004surface}, varifolds \cite{charon2013varifold,kaltenmark2017general} or normal cycles \cite{roussillon2019representation}. 

Assume that $S_0$ and $S_1$ are triangulated surfaces and that the cost function $D$ is replaced by a discrete approximation, still denoted $D$.
Then, the optimization problem can be reduced to one tracking explicitly the evolution of the vertices of the triangulation, using the reproducing kernel of $V$ denoted as $K$. This kernel is a matrix-valued function of two variables $x,y\in \mR^3$ such that, for all $\alpha, y\in \mR^3$, the vector field $x\mapsto K(x,y) \alpha$ belongs to $V$ and for all $v\in V$,
\[
\scp{v}{K(\cdot, y)\alpha}_V = \alpha^T v(y)
\]
where the left-hand side denotes the inner product in $V$.

Denote as $q_0 = (q_0(1), \allowbreak \ldots, q_0(N))$ the vertices of $S_0$. The reduced problem is expressed in terms of  evolving vertices $q(t) = (q(t,1), \ldots, q(t,N))$ and vectors $\alpha(t) = (\al(t,1), \ldots, \al(t,N))$, $t\in [0,1]$, minimizing (letting $S(t)$ denote the triangulated surface with vertices $q(t)$ and same topology (faces) as $S_0$)
\[
\int_0^1 \sum_{k,l=1}^N \alpha(t,k)^T K(q(t,k), q(t,l))\alpha(t,l) \, dt + D(S(1), S_1)
\]
subject to $q(0) = q_0$ and
\[
\partial_t q(t,k) = \sum_{l=1}^N K(q(t,k), q(t,l))\alpha(t,l)
\] for $k=1, \ldots, N$.
Moreover, the optimal vector field at time $t$ is given by
\begin{equation}
\label{eq:v}
v(t, \cdot) = \sum_{l=1}^N K(\cdot, q(t,l))\alpha(t,l).
\end{equation}

The interest of this formulation is that the trajectories $t\mapsto q(t, k)$ for $k=1, \ldots, N$ directly provide the streamlines starting from the vertices $q_0(k)$, $k=1, \ldots, N$ of the triangulation of $S_0$.

\begin{figure}
\centering
\includegraphics[trim=10cm 5cm 10cm 5cm,clip,width=0.6\textwidth]{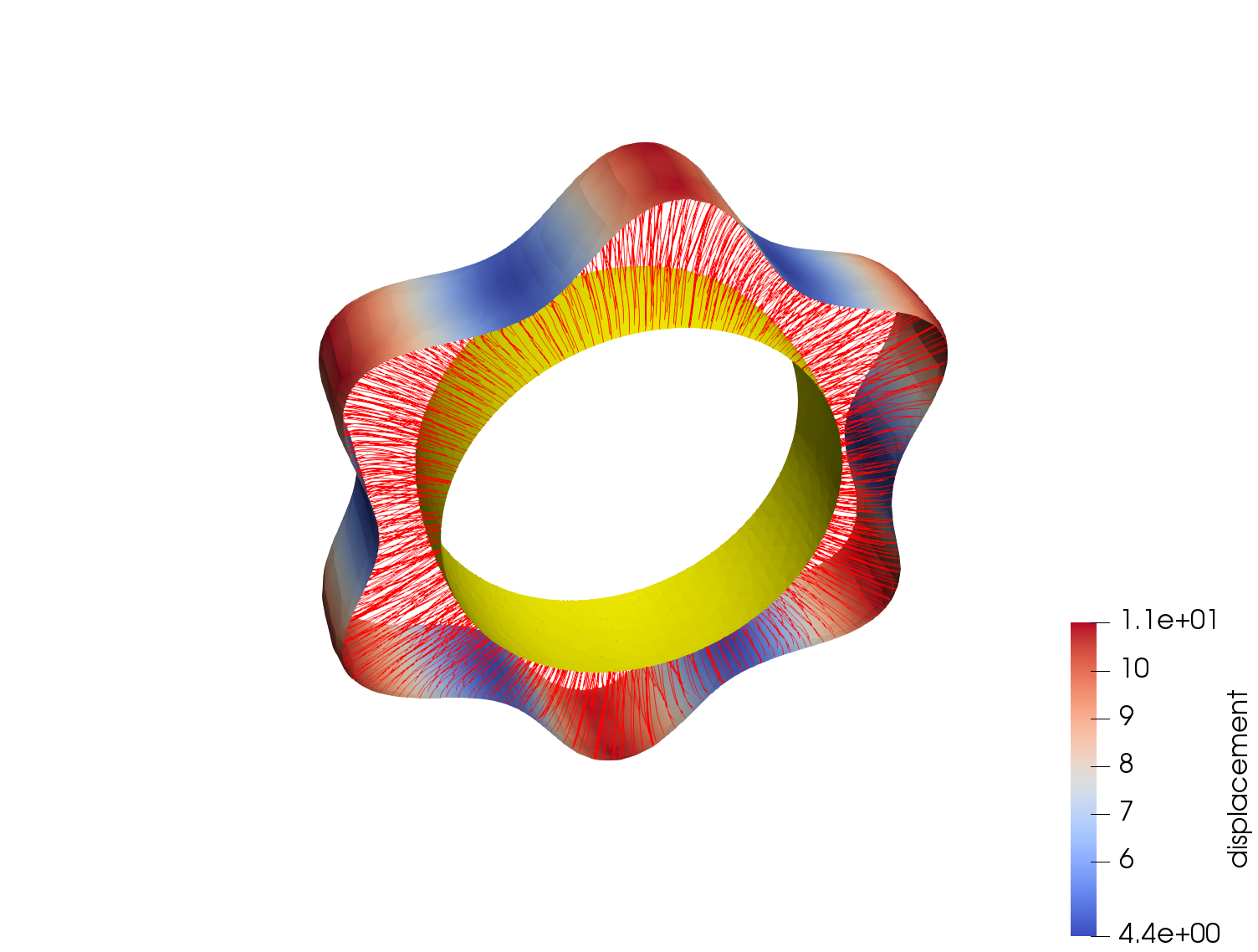}
\caption{\label{fig:flower} Synthetic data for thickness estimation between an inner ring of fixed radius and an outer one of variable radius.
%provided by the streamline lengths.
The colors on the outer ring show thickness values increasing from blue to red.}
\end{figure}

\medskip
This algorithm is modified by imposing a constraint ensuring that these streamlines are perpendicular to the evolving layers \cite{ratnanather20193D}. In the continuous setting, where, for each $t\in [0,1]$ $q(t, \cdot)$ is defined on the manifold $M$, this constraint can be formulated as $v(t, q(t,s)) = \lambda(t,s) \nu_{S(t)}(q(t,s))$, $t\in[0,1]$, $s\in M$. Here $\nu_S(x)$ denotes the (positively oriented) unit normal to an oriented surface  $S$ at $x\in M$, and $\la(t,s)$ is a scalar that we assume to be non-negative (with a proper orientation of $S(t)$ to prevent the trajectories from backtracking). The constraint is implemented in the equivalent form
\[
\sqrt{v(t,q(t,s))^Tv(t,q(t,s))} - \nu_{S(t)}(q(t,s))^Tv(t,q(t,s)) = 0,
\]
that is  discretized on the (evolving) triangulation of $S(t)$. The resulting constrained optimization problem is solved using an augmented Lagrangian method \cite{nocedal1999numerical} with each gradient descent step implemented using a limited-memory BFGS method (for details of methods using LDDMM with constraints, see \cite{arguillere2015shape}). Figure \ref{fig:flower} illustrates a synthetic example of how the thickness map can be estimated using this method.

Note that the LDDMM algorithm implicitly provides a flow of diffeomorphisms on the whole space $\mR^3$, given by solution of $\prt_t \phi(t, \cdot) = v(t, \phi(t, \cdot))$.
\medskip

Even if we will not use it in the following, it is interesting to notice is that this method generally provides a level-set formulation of the laminar coordinates. Indeed, under the mild assumption that $v(t,q(t,s))$ never vanishes, the function $\psi: [0,1]\times S_0\to \mR^3$ defined by $\psi(t, q_0(s)) = q(t,s)$ is an immersion. In practice, this mapping is in addition one-to-one and one can define without ambiguity a scalar function $F$ on a domain sandwiched by the surfaces $S_0$ and $S_1$ (or, more precisely, by $S_0$ and $S(1) \simeq S_1$),by
\[
F(q(t,s)) = t
\]
for $s\in M$. By construction, the streamlines are perpendicular to the level set of this function.
\medskip

\section{Equivolumetric Coordinates}
\label{sec:bok}
\subsection{Strict Condition}
As mentioned above, Bok's hypothesis requires that the cortical layers satisfy an equivolumetric constraint. Define a ``cortical tube'' as the volume delimited  by the cortical columns stemming from  the inner (grey-white) matter surface to the outer (pial) surface.
The hypothesis requires that the volume delimited by the intersection of cortical tubes and cortical surfaces remains roughly constant when the base patch is ``translated'' along the inner surface. We want to formalize this into a definition of  ``equivolumetric laminar coordinates.''

\begin{figure}
    \centering
    \includegraphics[trim=1cm 1cm 3cm 1cm, clip,width=0.95\textwidth]{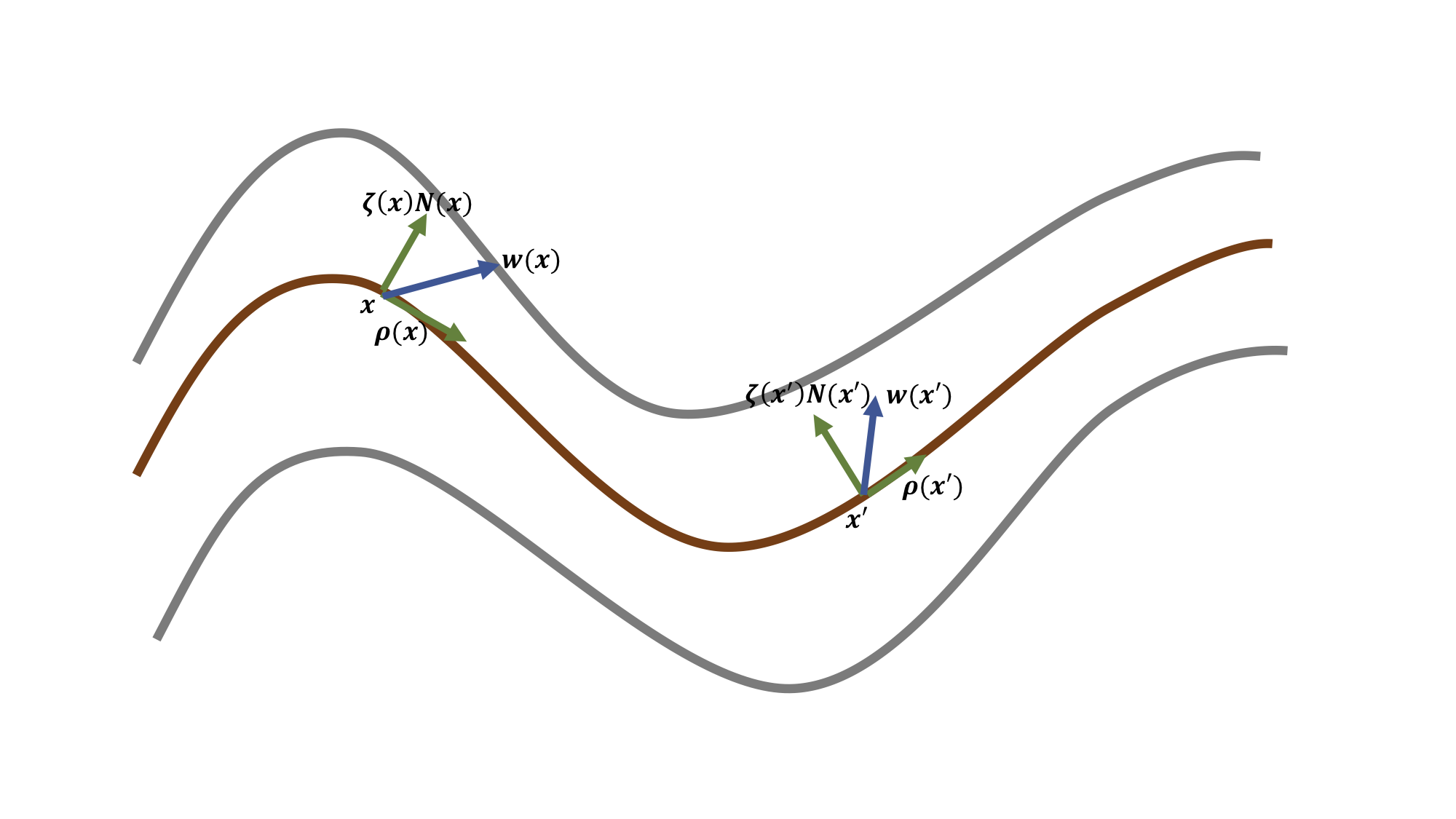}
    \caption{Illustration of the notation used in Proposition \ref{prop:sigma}. The evolution of intermediate surfaces (red) is captured by a vector field $w$ that decomposes into a tangential part $\rho$ and a normal part $\zeta N$.}
    \label{fig:prop1}
\end{figure}

Consider a laminar coordinate system $\psi: [0,1] \times S_0 \to \mR^3$. For $x\in S_0$, consider a small surface element $\delta S_0$ located at $x$ and the infinitesimal tube $\psi([0,1]\times S_0)$. Introduce a local chart $m: U \to \delta S_0$ on $\delta S_0$, where $U$ is an open subset of $\mR^2$. Let $\psi_m(t, \alpha, \beta) = \psi(t, m(\alpha, \beta))$. Then the  volume of the tube between layers $t_0$ and $t_1$ is given by
\begin{multline*}
\int_{t_0}^{t_1} \int_{U} \mathrm{det}(\prt_t \psi_m, \prt_\alpha \psi_m, \prt_\beta \psi_m)\, d\alpha\, d\beta\, dt \\ \simeq \int_{t_0}^{t_1} \int_{\delta S_0} \prt_t\psi(t,x)^T \nu(t,x) \sigma(t,x) d\mathrm{vol}_{S_0}(x)dt.
\end{multline*}
Here, we have used the notation $\nu(t,x)= \nu_{S(t)}(\psi(t,x))$. We also have denoted by $d\mathrm{vol}_{S_0}$ the volume form on $S_0$, which is given, in the local chart, by $|\prt_\al m \times \prt_\beta m|d\alpha d\beta$, and by $\sigma(t,x)$ the surface Jacobian induced by $\psi(t, \cdot): S_0 \to S(t)$ (infinitesimal ratio of area), defined in the chart by
\[
\frac{|\prt_\al \psi_m \times \prt_\beta \psi_m|}{|\prt_\al m \times \prt_\beta m|}\,.
\]
One has the following result that describes the evolution of $\sigma$ as a function of $t$.
\begin{proposition}
\label{prop:sigma}
Let $w: \Omega \mapsto \mR^3$ be defined by $w(\psi(t,x)) = \prt_t \psi(t,x)$ for all $(t,x) \in [0,1]\times S_0$. Define $N: \Omega \mapsto \mR^3$ by $N(\psi(t,x)) = \nu(t,x)$ and
decompose $w$ in the form
\[
w(y) = \rho(y) + \zeta(y) N(y),
\]
with $\rho(y) \perp N(y)$ for all $y\in \Omega$. Let $\rho_t$ denote the restriction of $\rho$ to $S(t)$.
One has
\begin{equation}
\label{eq:sigma}
\sigma^{-1} \prt_t \sigma = (\mathrm{div}_{S(t)} \rho_t - 2\zeta H_{S(t)}) \circ \psi
\end{equation}
where $\mathrm{div}_{S(t)}$ is the divergence operator on $S(t)$ and $H_{S(t)}$ the mean curvature on the same surface.
\end{proposition}
The notation used in this proposition is illustrated in Fig. \ref{fig:prop1}.
Recall that the divergence operator of a vector field $\rho$ on a surface $S$ is (for $p\in S$)
\begin{equation}
\label{eq:div.diff}
\mathrm{div}_S(\rho)(p) = e_1^T D\rho(p)e_1 + e_2^T D\rho(p)e_2
\end{equation}
where $e_1, e_2$ is any orthonormal basis of $T_pS$ (the tangent space to $S$ at $p$). Using the same basis, the mean curvature is given by
\begin{equation}
\label{eq:H.diff}
- 2 H_S(p) = e_1^T D\nu_S(p)e_1 + e_2^T D\nu_S(p)e_2.
\end{equation}
When convenient,  will use the notation: $H(t,x) = H_{S(t)}(\psi(t,x))$ to represent the mean curvature along a streamline.
\begin{proof}
We make the computation in a chart $m: U\subset \mR^2 \to S_0$ and let
\[
J(t,\alpha, \beta) = |\prt_\alpha \psi_m \times \prt_\beta \psi_m|
\]
so that $\sigma(t, m(\alpha, \beta)) = J(t,\alpha, \beta) / J(0, \alpha, \beta)$. Then
\[
\prt_t J = (\prt_t \prt_\alpha \psi_m \times \prt_\beta \psi_m)^T \nu_m + (\prt_\alpha \psi_m \times \prt_t \prt_\beta \psi_m)^T \nu_m
\]
with
\[
\nu_m(t, \alpha, \beta) = \nu(t, m(\alpha, \beta)) = \frac{\prt_\alpha \psi_m \times \prt_\beta \psi_m}{|\prt_\alpha \psi_m \times \prt_\beta \psi_m|}.
\]
Let $e_\alpha = \prt_\alpha \psi_m$, $e_\beta = \prt_\beta \psi_m$. Note that, letting $w_m(t, \alpha, \beta) = w(\psi_m(t, \alpha, \beta))$, one has
\[
w_m(t, \alpha, \beta) = \prt_t \psi_m(t, \alpha, \beta),
\]
so that
\[
\prt_t J = (\prt_\alpha w_m \times e_\beta + e_\alpha\times \prt_\beta w_m)^T\nu_m.
\]
Write $w_m = \rho_m + \zeta_m \nu_m$ so that
\[
\prt_t J = (\prt_\alpha \rho_m \times e_\beta + e_\alpha\times \prt_\beta \rho_m)^T\nu_m + \zeta_m (\prt_\alpha \nu_m \times e_\beta + e_\alpha\times \prt_\beta \nu_m)^T\nu_m.
\]
One has
\[
\nu_m^T(\prt_\alpha \rho_m \times e_\beta + e_\alpha\times \prt_\beta \rho_m) = J (\text{div}_{S(t)} \rho_t)(\psi_m(t, \cdot))
\]
and
\[
\nu_m^T(\prt_\alpha \nu_m \times e_\beta + e_\alpha\times \prt_\beta \nu_m) = -2J H_{S(t)}(\psi_m(t, \cdot)).
\]
These identities can be proved from \eqref{eq:div.diff} and \eqref{eq:H.diff} by introducing an orthonormal basis $(e_1,e_2)$ of $T_{\psi_m(t,x)}S(t)$ and expanding the left-hand sides as functions of the coordinates of $e_\al$ and $e_\beta$ in this basis (cf. \cite{younes2019shapes}, lemma 3.19). Using this, we get
\[
\prt_t J = J\, (\mathrm{div}_{S(t)}\rho_t - 2H_{S(t)}) \circ \psi_m
\]
from which one deduces that
\[
\prt_t \sigma = \sigma\,(\mathrm{div}_{S(t)}\rho_t - 2H_{S(t)}) \circ \psi
\]
\end{proof}
\bigskip

Define the ``equivolumetric thickness'' along the streamline starting at $x$ by
\begin{equation}
    \label{eq:equiv.thick}
\gamma(t,x) = \int_0^t \prt_t\psi(u,x)^T \nu(u,x) \sigma(u,x)\, du.
\end{equation}

A strict interpretation of Bok's hypothesis requires that this expression does not depend on $x$, i.e., that there exists a function $t\mapsto \lambda(t)$ such that
\begin{equation}
\label{eq:bok0}
\prt_t\psi(t,x)^T \nu(t,x) \sigma(t,x) = \lambda(t)
\end{equation}
for all $x\in S_0$. (In which case $\gamma(t,x) = \int_0^t \lambda(u) du$.)
One can, without loss of generality, assume that $\lambda$ does not depend on $t$, which can be achieved by applying a time change to the evolution. More precisely, one can replace $\psi$ by $\tilde\psi$ such that
$\tilde \psi(\tau(t), x) = \psi(t, x)$
with
\[
\tau(t) = \frac{\int_0^t \lambda(u) du}{\int_0^1 \lambda(u) du}
\]
in which case $\tilde \psi$ satisfies \eqref{eq:bok0} with constant right-hand side
\[
\tilde \lambda = \int_0^1 \lambda(u) du.
\]

%We point out that the evolution of the surface $S(t)$ is characterized by the normal component of $\prt_t \psi$ (the tangential component only affect the induced parametrization).
Therefore assuming that $\lambda$ is constant, we now apply Proposition \ref{prop:sigma} to obtain a surface propagation equation that is equivalent to \eqref{eq:bok0}. We will decompose $w = \prt_t \psi \circ \psi^{-1}$ in the form
\[
w(y) = \rho(y) + \frac{\lambda}{\sigma(\psi^{-1}(y))} N(y),
\]
as required by \eqref{eq:bok0}. One then has the following proposition.
\begin{proposition}
\label{prop:bok1}
A laminar coordinate system $\psi$ satisfies Bok's hypothesis if there exists a vector field $\rho$ tangent to the layers defined by $\psi$ and a constant $\lambda$ such that
 \begin{equation}
    \label{eq:bok1}
    \left\{
    \begin{aligned}
    \prt_t \psi(t,x)  &= \rho(\psi(t,x))   + \frac{\lambda}{\sigma(t,x)} \nu(t,x) \\
\prt_t \sigma(t,x) &= \sigma(t,x) \mathrm{div}_{S(t)} \rho(\psi(t,x)) - 2\lambda H(t,x)
\end{aligned}
\right.
\end{equation}
with $\psi(0, x) = x$ and $\sigma(0, x) = 1$ for all $x\in S_0$.
\end{proposition}
Equation \eqref{eq:bok1}
provides an evolution equation controlled by  $\rho$ (such that at all times, $\rho(t, \cdot)$ is a vector field on the evolving surface $S(t)$) whose solution satisfies the equivolumetric hypothesis. A detailed study of this system of equations (including its well-posedness for a given choice of $\rho$), and of the optimal control problem consisting in optimizing $\rho$ with the constraint that $S(1) = S_1$ are challenging open problems that will not be addressed in this paper. Even if possible, it would also be counter-intuitive to require a constant equivolumetric thickness, $\gamma(1,x)$, in \eqref{eq:equiv.thick}. So in the next section, we discuss a solution to a simpler problem that we refer to as a ``localised'' Bok's hypothesis.

\subsection{Localised Bok's hypothesis}
In this section, we replace the strong constraint \eqref{eq:bok0} by a weaker one
\begin{equation}
\label{eq:bok2}
\prt_t\psi(t,x)^T \nu_{S(t)}(\psi(t,x)) \sigma(t,x) = \lambda(t) c_0(x)
\end{equation}
for a given function $c_0$. Here again, there is no loss of generality in assuming that $\lambda$ is constant, and, including if needed this constant in $c_0$, in taking $\lambda = 1$ (so that $c_0$ coincides with the equivolumetric thickness $\gamma(1,x)$). Proposition \ref{prop:bok1} can be directly extended in this setting.
\begin{proposition}
\label{prop:bok2}
A laminar coordinate system $\psi$ satisfies the localised form of Bok's hypothesis for a given equivolumetric thickness $c_0$ if there exists a vector field $\rho$ tangent to the layers defined by $\psi$ such that
 \begin{equation}
    \label{eq:bok3}
    \left\{
    \begin{aligned}
    \prt_t \psi(t,x)  &= \rho(\psi(t,x))   + \frac{c_0(x)}{\sigma(t,x)} \nu(t,x) \\
\prt_t \sigma(t,x) &= \sigma(t,x) \mathrm{div}_{S(t)} \rho(\psi(t,x)) - 2c_0(x) H(t,x)
\end{aligned}
\right.
\end{equation}
with $\psi(0, x) = x$ and $\sigma(0, x) = 1$ for all $x\in S_0$.
\end{proposition}

Because the function $c_0$ can be chosen freely, \eqref{eq:bok2} (or a solution of system \eqref{eq:bok3}) is much easier to obtain while satisfying the condition $\psi(1,S_0) = S_1$ and we now show that it can be achieved starting from any  laminar coordinate system $\psi$ by applying a space-dependent time change. Indeed, let $\tilde \psi(\tau(t,x),x) = \psi(t, x)$, where, for each $x\in S_0$, $t\mapsto \tau(t,x)$ is an increasing differentiable function from $[0,1]$ onto $[0,1]$. Then, introducing as above local coordinates $\alpha$ and $\beta$ on $S_0$, we have
\begin{align*}
\partial_\alpha \psi(t,x) &= \partial_\alpha \tau\, \prt_t \tilde \psi (\tau, x) + \prt_\alpha \tilde \psi(\tau, x)\\
\partial_\beta \psi(t,x)  &= \partial_\beta \tau\, \prt_t \tilde \psi(\tau,x) + \prt_\beta \tilde \psi(\tau, x)\\
\partial_t \psi(t,x) &= \partial_t\tau\, \prt_t\tilde \psi(\tau, x)
\end{align*}

As a consequence:
\begin{align*}
\nu(t,x)^T\prt_t \psi(t,x) \sigma(t,x)
&=\mathrm{det}(\prt_t \psi, \partial_\alpha \psi, \partial_\beta \psi)(t,x) \\
&= \partial_t\tau(t,x)\, \mathrm{det}(\partial_t \tilde \psi, \partial_\alpha \tilde \psi, \partial_\beta \tilde \psi)(\tau,  x)\\
&= \partial_t\tau(t,x)\, \tilde \nu(\tau, x)^T \prt_t\tilde \psi (\tau,x) \tilde \sigma(\tau,x).
\end{align*}
In order that \eqref{eq:bok2} holds for $\tilde \psi$, we therefore need to define $\tau$ so that (for some function $c_0$)
\[
c_0(x) \partial_t \tau(t,x) =  \nu(t, \psi(t, x))^T \prt_t\psi(t, x) \sigma(t,x)
\]
In order to have $\tau(0,m) = 0$, $\tau(1, m) = 1$, we need
\[
c_0(x) = \int_0^1 \nu(u, \psi(u,x))^T \prt_t\psi(u, x) \sigma(u,x)\, du
\]
and then
\begin{equation}
\tau(t,x)  = \frac{1}{c_0(x)} \int_0^t \nu(u, \psi(u,x))^T \prt_t\psi(u, x) \sigma(u,x)\, du.
\end{equation}
\bigskip

So, the time change is provided by the relative volumetric depth the streamlines (that are left unchanged in the operation). The equivolumetric layers at level $\epsilon$ are provided by points $\psi(t,x)$ along the streamlines satisfying $\tau(t,x) = \epsilon$ (in the original parametrization).

\subsection{Interpretion of recent models}
We now interpret two recent attempts to model Bok's hypothesis \cite{waehnert2014anatomically,leprince2015combined} in our framework.

Waehnert et al. \cite{waehnert2014anatomically} start with streamlines estimated using \eqref{eq:waehnert1}. Equivolumetric layers are then estimated using the same equation,  replacing the constant value $\rho$ by a function $\rho(x)$, that is determined as follows. Using our notation, one first makes the assumption that
\[
\sigma(t,x) = (1-t) + t \sigma(1,x)
\]
therefore making a linear approximation of the surface change. Here, one takes $t = s/\theta$, where $s$ is the arc length along the streamline and $\theta$ is the thickness (the length of the streamline). For $x\in S_0$, the value of $\sigma(1,x)$ is estimated as a function of the curvatures at both ends of the streamline starting at $x$ (we refer to \cite{waehnert2014anatomically,Kemper2018} for more details and justification). Integrating along streamlines, which are perpendicular to the layers because of the level set formulation, one obtains an expression of the equivolumetric depth given by
\[
\mathcal V_x(\rho) = \theta \int_0^\rho \sigma(t,x) dt = \theta \rho\left( 1 + \frac\rho2 (\sigma(1,x)-1)\right).
\]
Given $\epsilon\in[0,1]$, one defines a target level $\rho_\epsilon(x)$ corresponding to the layer at equivolume $\epsilon$ by solving $\mathcal{V}_x(\rho) = \alpha \mathcal V_x(1)$, which is a quadratic equation in $\rho$.

\medskip

Leprince et al. \cite{leprince2015combined} start with a Laplacian-based level-set definition of streamlines \cite{jones2000three}. They estimate $\sigma(t, \cdot)$ along the streamlines by solving
\[
\prt_t \sigma(t, x) = - 2 \theta(x) \sigma(t,x) H(t,x)
\]
which corresponds to \eqref{eq:sigma} for a constant-speed normal evolution $\prt_t \psi(t,x) = \theta(x) \nu(t,x)$. Note that, in the level set approach, one has $\nu(t,x) = (\nabla F/|\nabla F|)\allowbreak(\psi(t,x))$ and $2H(t,x) = -\mathrm{div}(\nabla F/|\nabla F|)(\psi(t,x))$. Equivolumetric layers are then deduced from this computation.

\subsection{Numerical Implementation}
Because they rely on level sets, the two recent approaches \cite{waehnert2014anatomically,leprince2015combined} are Eulerian, i.e., they work in the 3D volume $\Omega$ and layers are isosurfaces associated with scalar functions defined on $\Omega$ while streamlines are integrated using the gradient of these functions.

Our approach is different in that it is directly modeling layers as parametrized surfaces $S(t) = \psi(t, S_0)$, so that our implementation is based on a triangulation of $S_0$ and a discretization of the time interval. The streamlines are obtained as solutions of the ODE $\prt_t y = v(t,y)$, where $v$ is obtained from the LDDMM algorithm and provided by \eqref{eq:v}. Importantly, this time-dependent vector field is discretized in time only, and known analytically as a function of $y$. In particular, its space derivatives can be evaluated without approximation. It also specifies a flow of diffeomorphisms of $\mR^3$ through the equation
\[
\prt_t \phi(t,x) = v(t, \phi(t,x))
\]
with $\phi(0,x) = x$ for all $x\in \mR^3$.

The surface Jacobian $\sigma$ can be evaluated using the evolving triangulated surfaces: if $x$ is a vertex on $S_0$, we let $a(0,x)$ denote the area of the one-ring centered at $x$ (the union of all triangles that contain $x$). Similarly, we let $a(t,x)$ denote the area of the one-ring around $\phi(t,x)$ in the triangulated surface $S(t) = \phi(T, S_0)$ (which has the same triangle structure as $S_0$). On can then define
\[
\sigma(t,x) = \frac{a(t,x)}{a(0,x)}.
\]
This is the approximation that are used in our simulations, and it is accurate provided that the triangulation of $S_0$ is fine enough without flat triangles. An alternative procedure is also possible, since one has, in this context
\begin{equation}
\label{eq:sigma.lddmm}
\sigma(t,x) = \mathrm{det}(\prt_x\phi(t,x)) |\prt_x\phi(t,x)^{-T} \nu(0,x)|
\end{equation}
for $x\in S_0$. (The ``$-T$'' exponent refers to the inverse of the transpose matrix.) Recall that $\phi(t, \cdot)$ is (for fixed time $t$) a diffeomorphism of $\mR^3$, hence defined on the whole space (unlike $\psi$, which, for laminar coordinates, is only defined on $S_0$, and in this special case, is defined as the restriction of $\phi$ to this surface). This implies that $\prt_x\phi(t,x)$ is a $3\times3$ matrix. To prove \eqref{eq:sigma.lddmm} one can just notice that, in a local chart $m$
\[
\sigma(t,x) = \frac{|(\prt_x \phi(t,x) \prt_\alpha m) \times (\prt_x \phi(t,x) \prt_\beta m)|}{|\prt_\alpha m \times \prt_\beta m|}
\]
and  use the fact that for any matrix $A$ and vector $u$ and $v$, one has $Au\times Av = \mathrm{det}(A) A^{-T}(u\times v)$.

The time evolution of the vector $\zeta(t,x) = \mathrm{det}(\prt_x\phi(t,x)) \prt_x\phi(t,x)^{-T} \nu(0,x)$ is provided by
\[
\prt_t \zeta =  \mathrm{div}(v) (\phi(t, x)) - \prt_x v(\phi(t, x)) \zeta(t,x)\,.
\]
This can be integrated along streamlines using the expression of $v$ in \eqref{eq:v}, from which, as mentioned, space derivatives can be evaluated exactly.

\section{Results}
\label{sec:results}
\begin{figure}
    \centering
\includegraphics[trim=5cm 0cm 5cm 1.5cm, clip, width=0.75\textwidth]{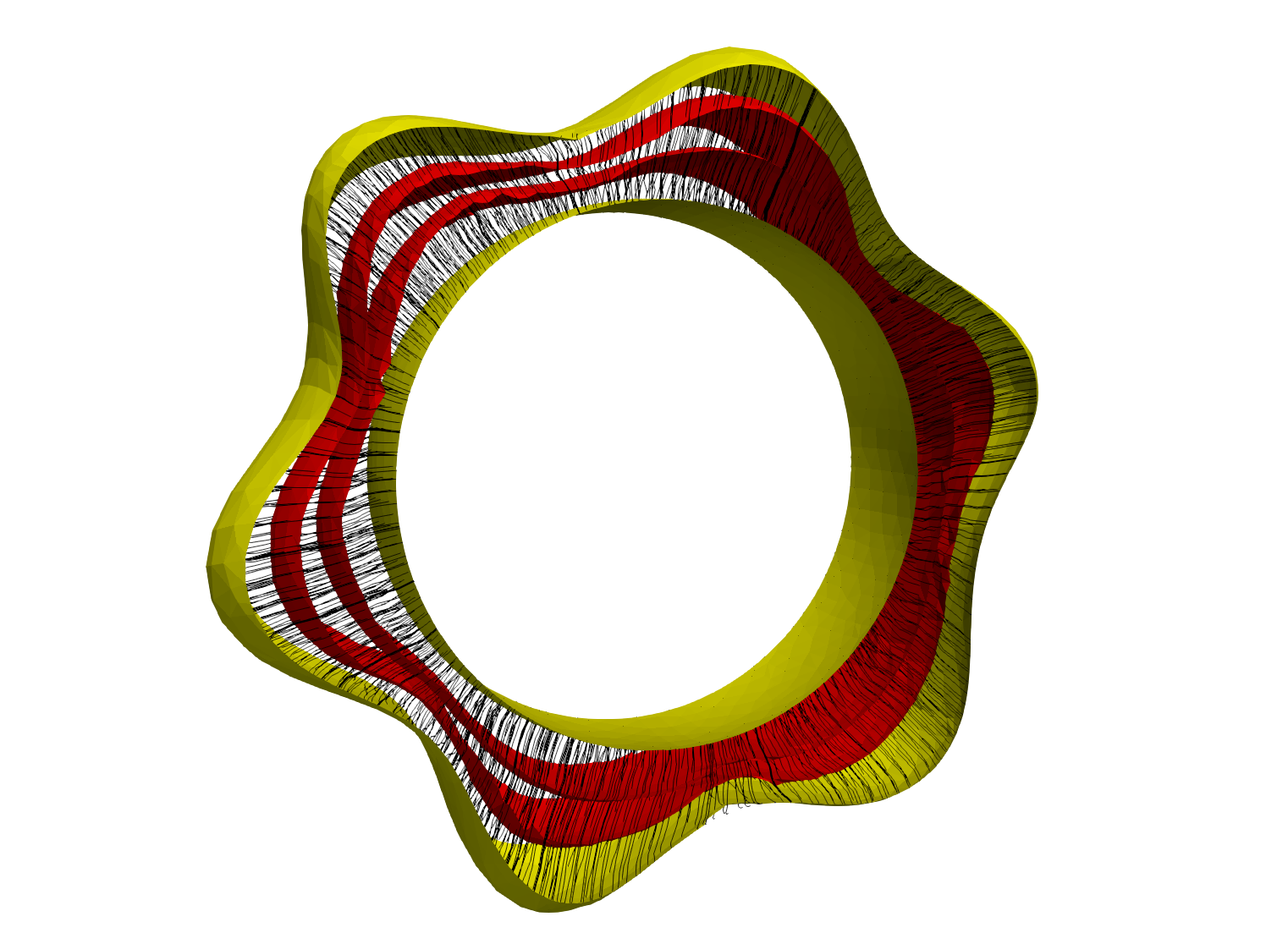}
    \caption{\label{fig:flower2} Estimated equivolumetric layers (in red) at times $t=0.3$ and $0.6$ for the synthetic example of Figure \ref{fig:flower}}.
\end{figure}

\begin{figure}
    \centering
\includegraphics[trim=0cm 15cm 0cm 0cm, clip,width=0.9\textwidth]{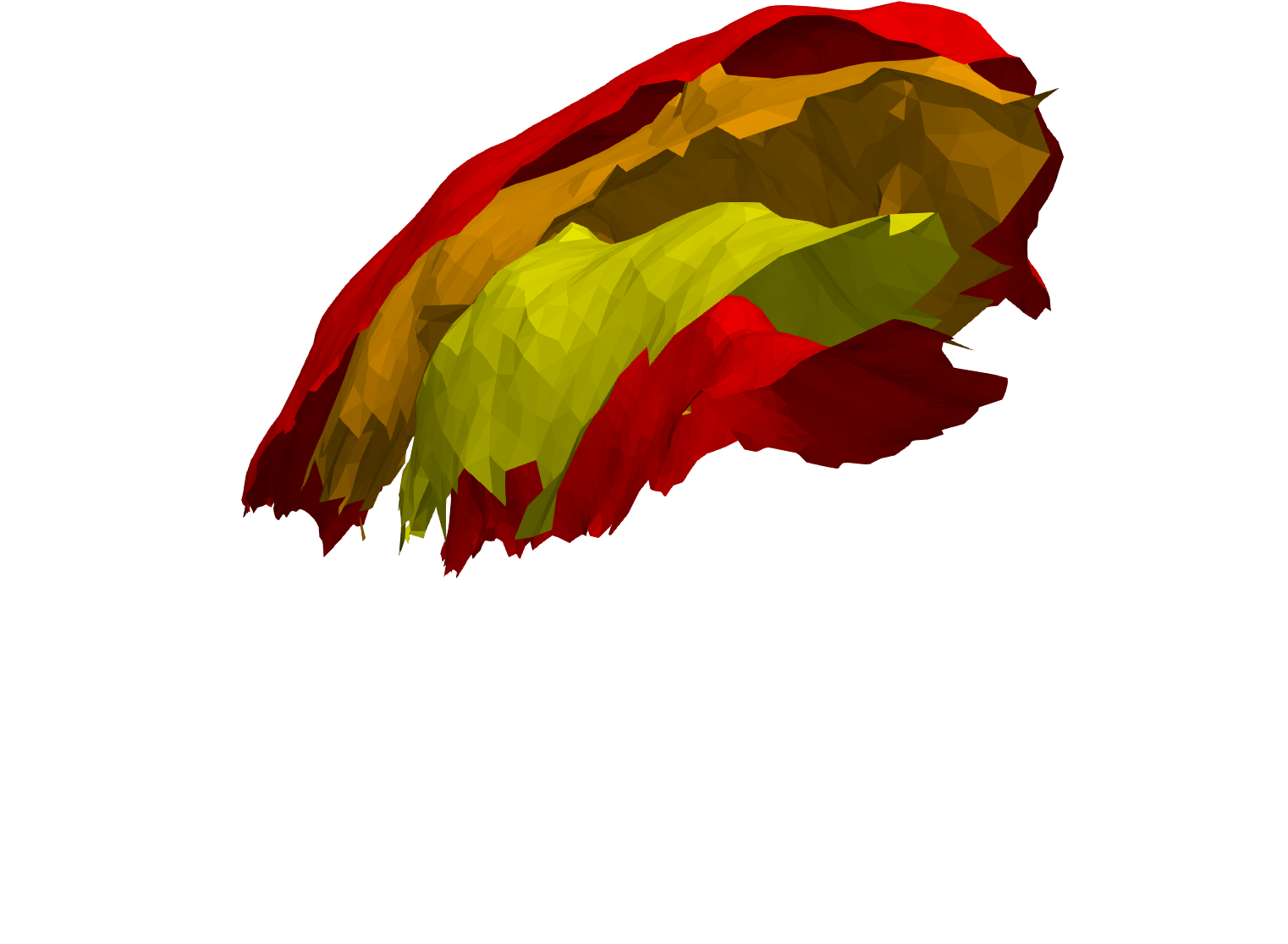}
\includegraphics[trim=0cm 15cm 0cm 0cm, clip,width=0.9\textwidth]{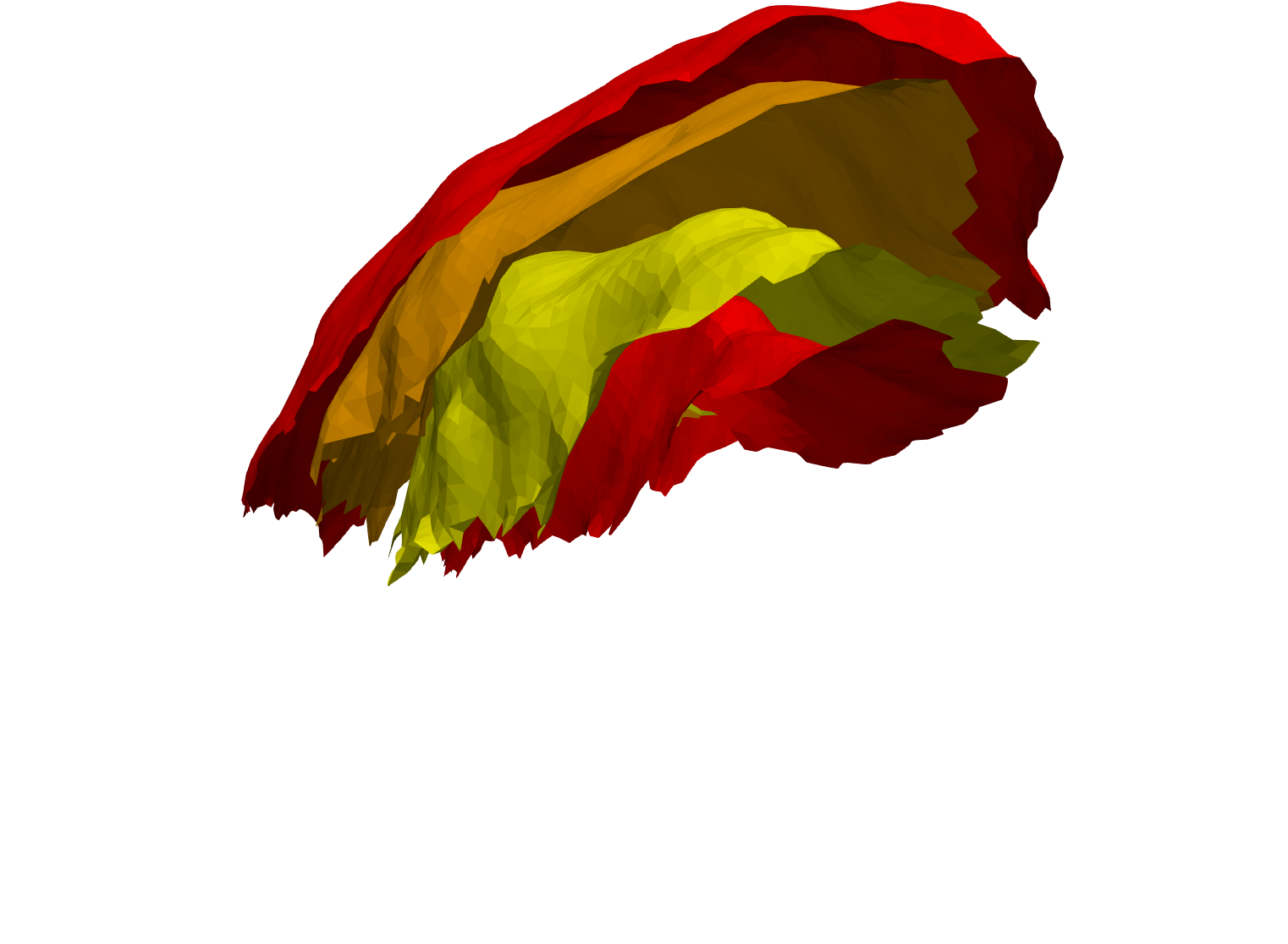}
\includegraphics[trim=0cm 15cm 0cm 0cm, clip,width=0.9\textwidth]{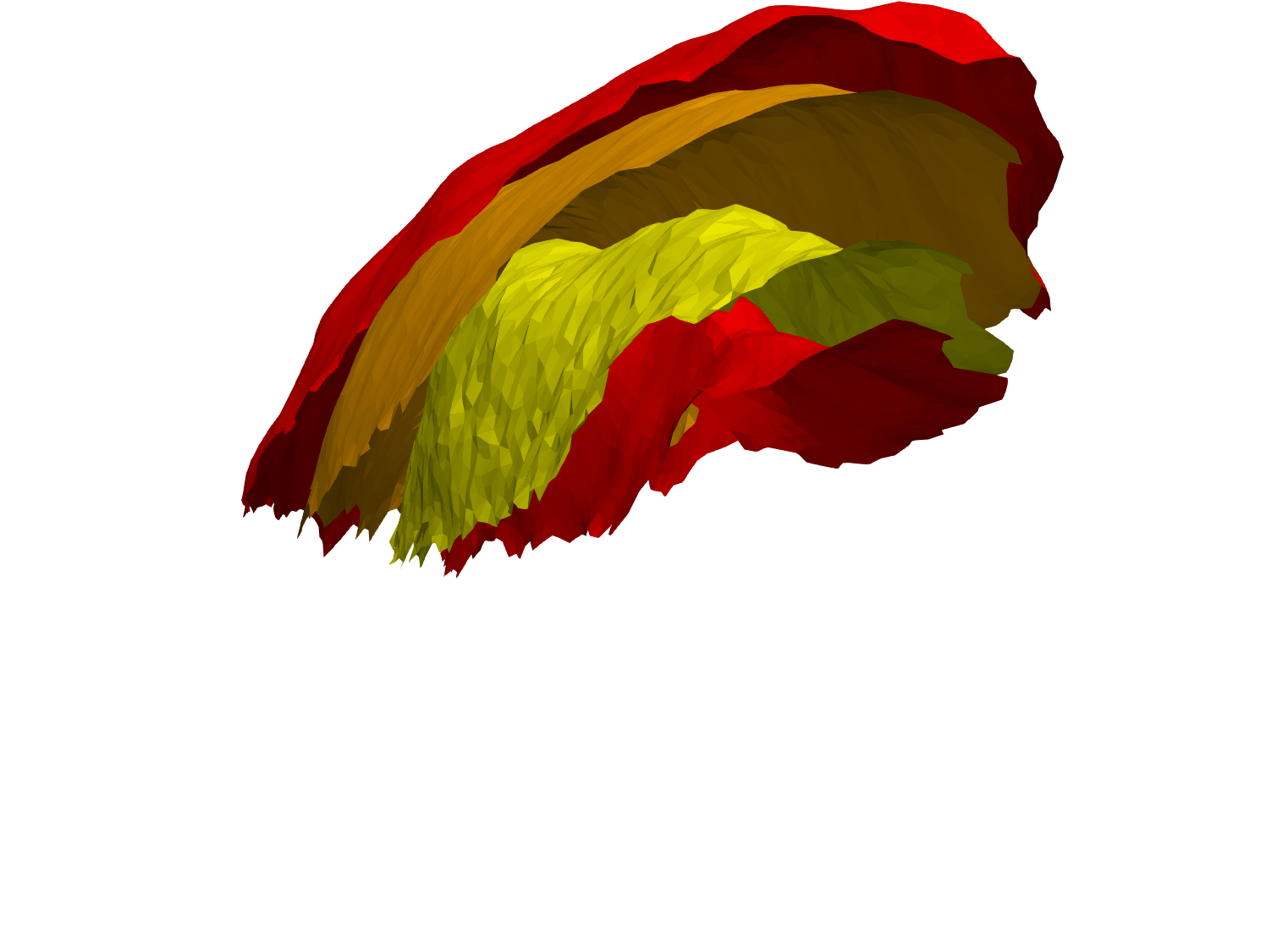}
    \caption{\label{fig:marmoset} Equivolumetric layers for a marmoset auditory cortex at $t=0.25$ (orange) and $t=0.75$ (yellow) using Laplacian \cite{leprince2015combined, leprince2017highres} (top), level set \cite{waehnert2014anatomically, hunterburg2017github} (middle) and the proposed method (bottom).}
\end{figure}

\begin{figure}
    \centering
\includegraphics[trim=0cm 10cm 0cm 5cm, clip,width=0.9\textwidth]{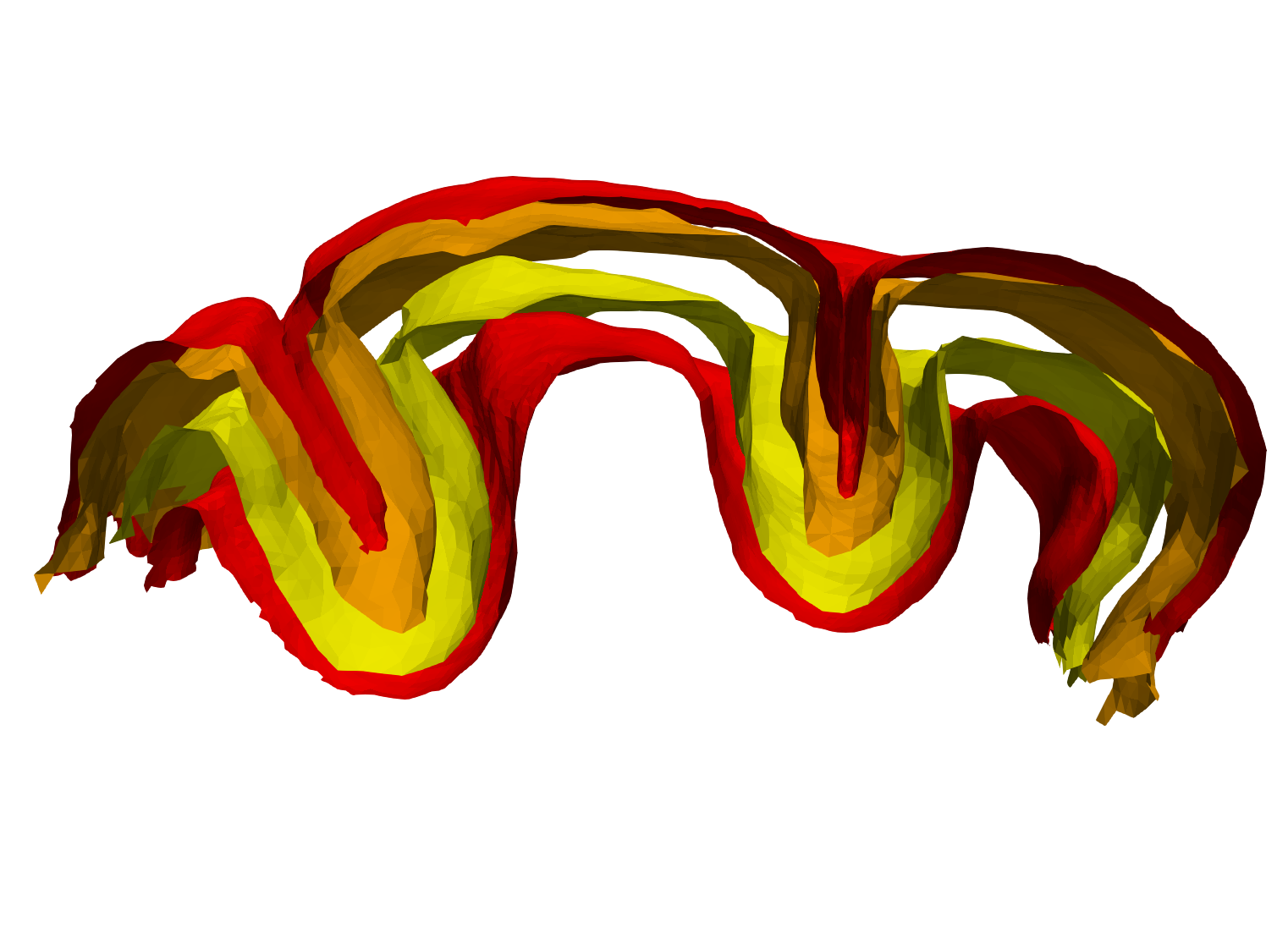}
\includegraphics[trim=0cm 10cm 0cm 5cm, clip,width=0.9\textwidth]{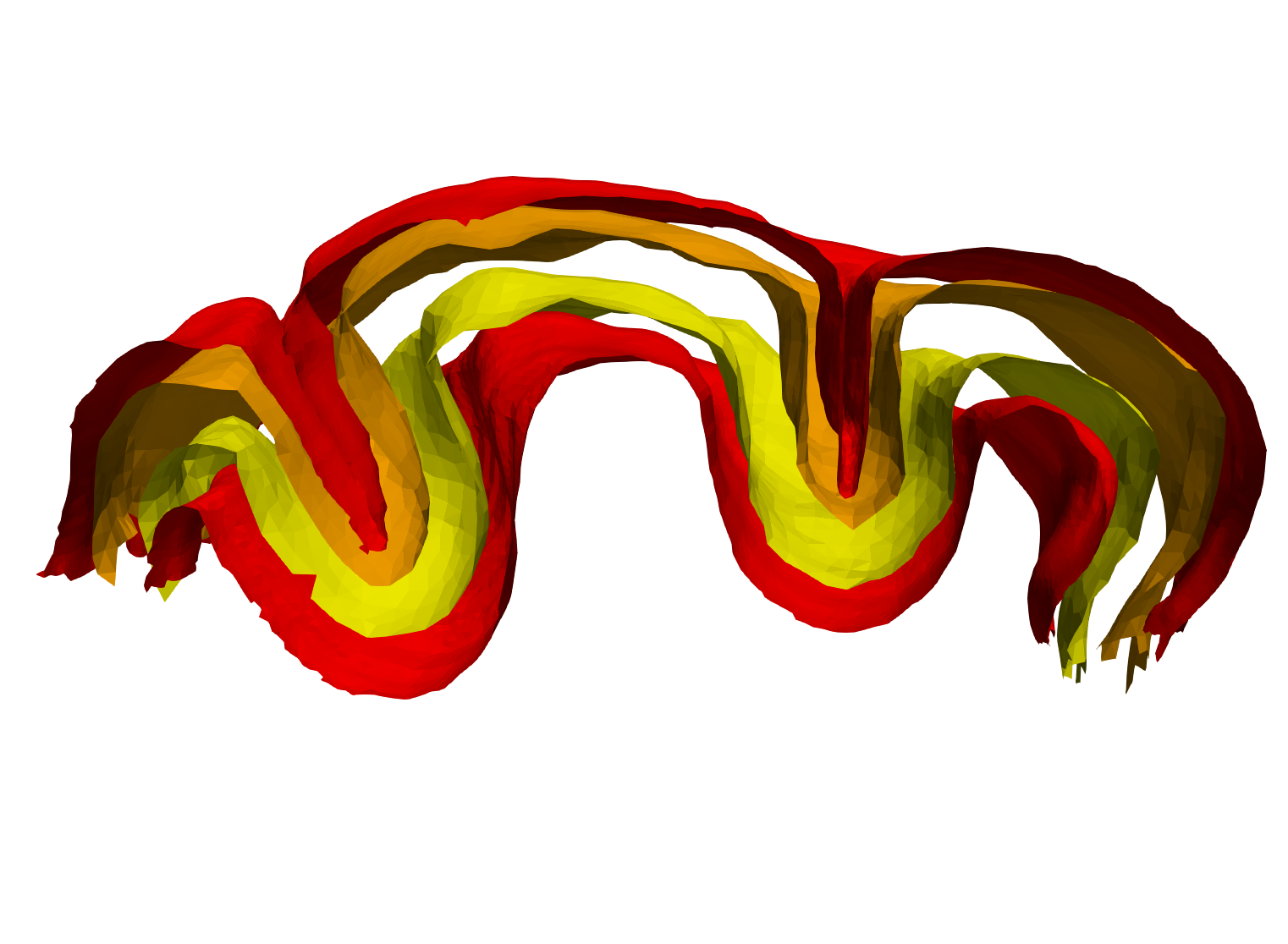}
\includegraphics[trim=0cm 10cm 0cm 5cm, clip,width=0.9\textwidth]{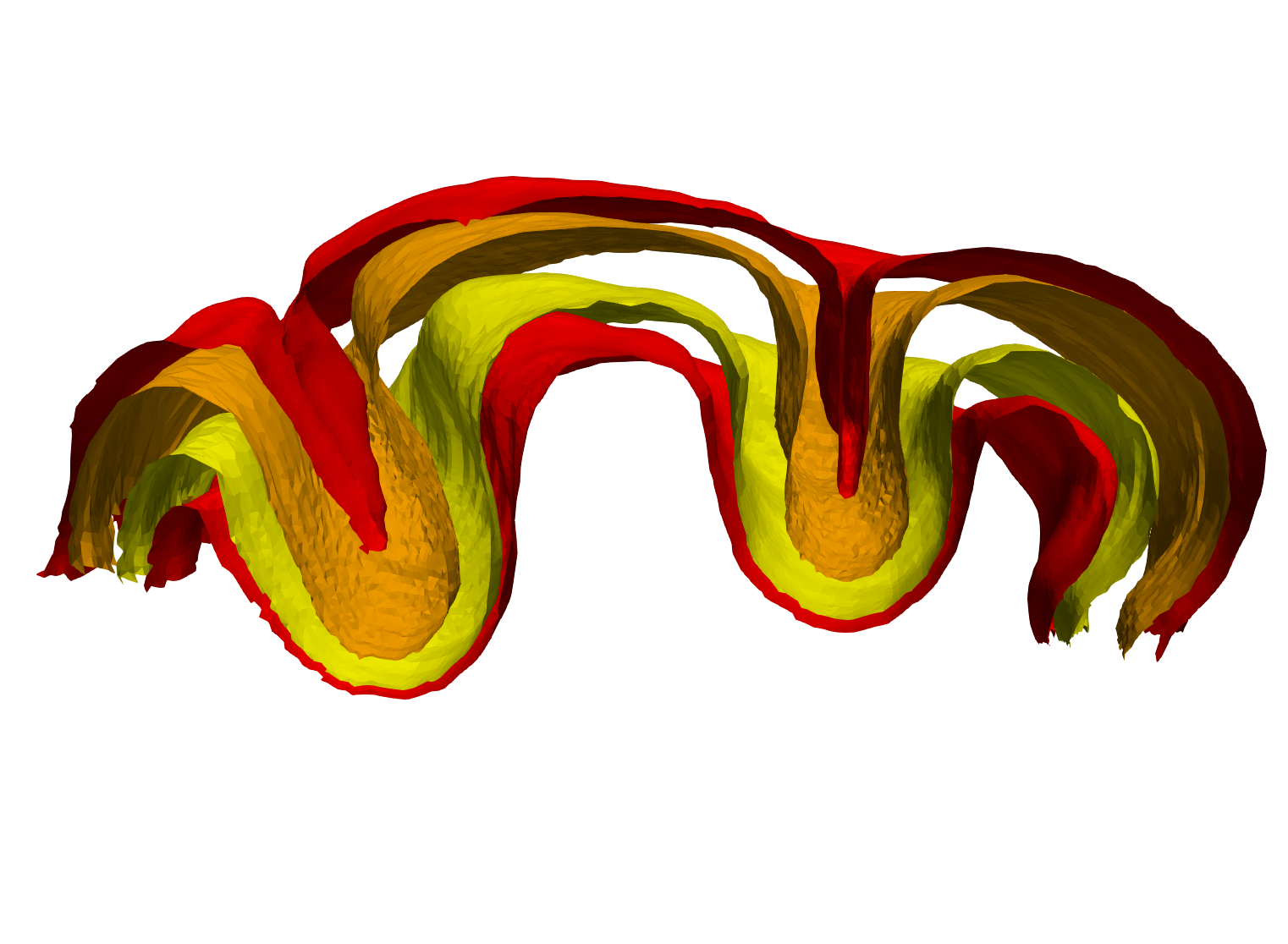}
    \caption{\label{fig:cat} Equivolumetric layers for a feline auditory cortex at $t=0.25$ (orange) and $t=0.75$ (yellow) using Laplacian \cite{leprince2015combined, leprince2017highres} (top), level set \cite{waehnert2014anatomically, hunterburg2017github} (middle) and the proposed method (bottom).}
\end{figure}

\begin{figure}
    \centering
    \begin{tabular}{cc}
\includegraphics[width=0.45\textwidth]{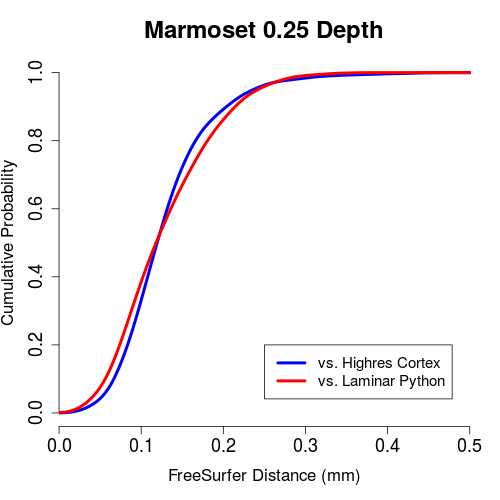} & \includegraphics[width=0.45\textwidth]{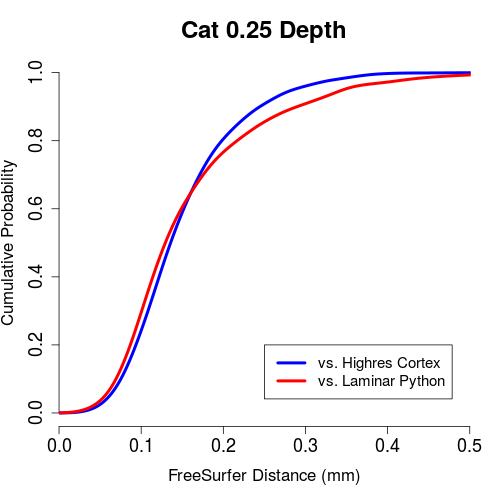} \\
\includegraphics[width=0.45\textwidth]{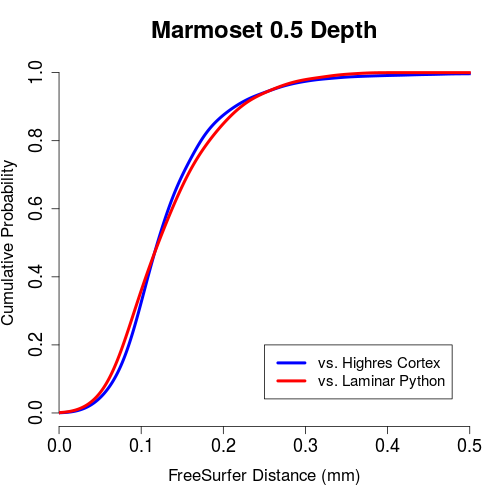} & \includegraphics[width=0.45\textwidth]{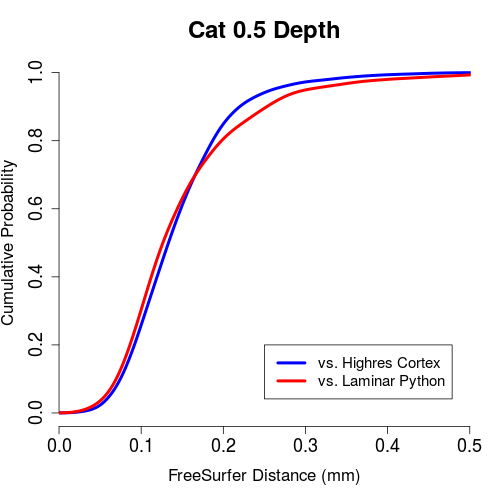} \\
\includegraphics[width=0.45\textwidth]{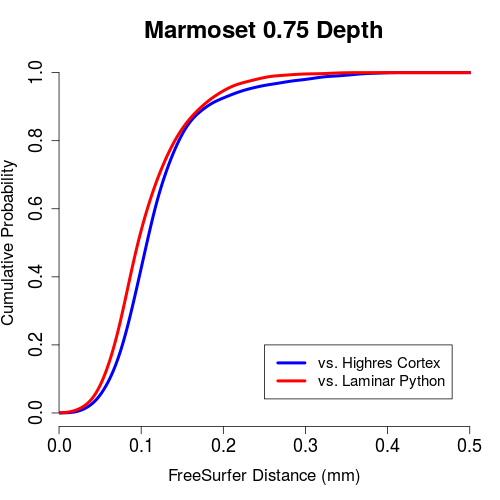} & \includegraphics[width=0.45\textwidth]{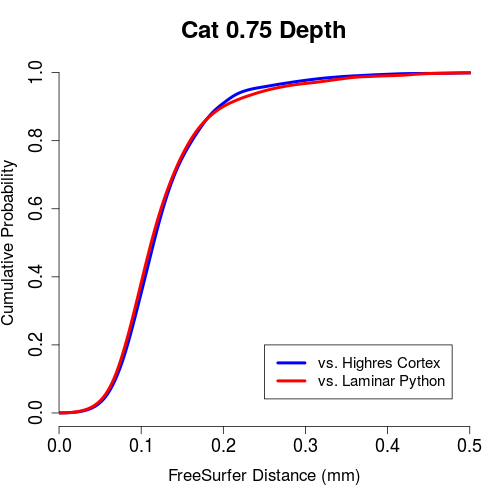}
\end{tabular}
    \caption{\label{fig:comp} CDFs of FreeSurfer distances of equivolumetric surfaces at $t=0.25, 0.5, 1.0$ from Figs. \ref{fig:marmoset} and \ref{fig:cat} computed via Github packages - Laplacian (Highres Cortex, \cite{leprince2017highres}) and Level Set (Laminar Python, \cite{hunterburg2017github}) - relative to surfaces computed via the proposed method.}
\end{figure}

Results of the numerical implementation are presented for three cases. Figure \ref{fig:flower2} shows the result for the synthetic data of Figure \ref{fig:flower}. Figures \ref{fig:marmoset} and \ref{fig:cat} respectively show the results for the marmoset auditory cortex (obtained from \cite{Woodward2018}) and feline auditory cortical regions (obtained from \cite{Berger2017}). Here, the proposed method is compared with those of Waehnert et al. \cite{waehnert2014anatomically} and Leprince et al. \cite{leprince2015combined} computed via Github packages \cite{leprince2017highres,hunterburg2017github}. Figure \ref{fig:comp} shows the corresponding cumulative distribution of distances of equivolumetric surfaces at $t=0.25, 0.5, 1.0$ relative to those via the proposed method. Here, for surfaces $S_1$ and $S_0$ with vertices $x_i \in S_1$ and $y_j \in S_0$, the distance at the $i$\textsuperscript{th} vertex of $S_1$ is \mbox{$d_i = \frac{1}{2} (d(x_i, y_n) + d(x_m, y_n))$} where $n = \underset{j}{\arg \min} \; d(x_i, y_j)$ and $m = \underset{k}{\arg \min} \; d(x_k, y_n)$. This FreeSurfer distance \cite{Fischl2000} returns a value for every vertex and making it more robust against outliers that arise with the Hausdorff distance.

\section{Discussion}
A unified theoretical framework for describing laminar coordinate systems for cortical regions has been developed. Herein, several algorithms including an interpretation of volumeteric-based ones are offered such that equivolumetric layers consistent with Bok's hypothesis can be computed. 

Granted that in the marmoset the auditory cortex is the only area that is folded with a gyral crown, the protrusions of the upper layers seen with the method of Leprince et al. \cite{leprince2015combined} may be attributed to the diverging normal vector fields as one approaches the outer surface. This divergence warrants finer discretization which may be computationally expensive. The primary and higher-order auditory cortical regions in the cat reveals greater differences between the three methods. The $t=0.25$ layer appears to be closer to the sulcal fundi for Leprince et al. and our methods; in contrast the $t=0.75$ layer appears to be closer to the gyral crowns for Waehnert et al. and our methods. These differences can be quantified via a distance metric (Figure \ref{fig:comp}) and may be attributed to the representation of curvature --direct or indirect-- in the computations. Future work will examine the how disease and disorder affect the equivolumetric depths of layers, pyramidal cells and other cortical elements in 3D. This will build upon previous work in 2D \cite{Berger2017}. Our laminar coordinate system also has a straightforward application to cortical layer segmentation. As a prior, it could help distinguish between layers of similar microanatomy
(e.g cortical layers III and V) that can not be separated using staining intensity alone.

The time change introduced to reparametrize the coordinates in order to make them compliant with Bok's hypothesis left the streamlines invariant while changing the layers. As a consequence, if streamlines were perpendicular to the layers to start with, this property is generally lost after reparametrization. Finding an equivolumetric coordinate system with perpendicular streamlines is a significantly more arduous problem. Using Proposition \ref{prop:bok2}, in which one must set $\rho=0$, one sees that this problem requires to estimate a scalar field $c_0$ on $S_0$ such that the solution of \eqref{eq:bok3} satisfies $\psi(1, S_0) = S_1$. Whether this inverse problem is well posed, and whether stable numerical algorithms can be designed to solve it, are open questions that we plan to address in future work.

\section*{Acknowledgements}
Support from National Institutes of Health (P41-EB015909, R01-DC016784) is gratefully appreciated. The feline data was obtained from Professor Kral's group via grant 01GQ1703 from Federal Ministry of Education and Research of Germany (BMBF).

\bibliographystyle{siam}
\bibliography{references}
\end{document}